\documentclass[12pt,twoside,reqno]{amsproc}

\usepackage[centertags]{amsmath}

\usepackage{amssymb,amsfonts,amsthm}
\newcommand{\TryPackage}[3]{\IfFileExists{#1.sty}
{\usepackage{#1}#2}{#3}}\TryPackage{mathrsfs}
{\renewcommand{\mathcal}{\mathscr}}

{
\TryPackage{eucal}{}{}}
\newtheorem{theorem}{Theorem}[section]
\newtheorem{cor}[theorem]{Corollary}
\newtheorem{lemma}[theorem]{Lemma}
\newtheorem{prop}[theorem]{Proposition}

\theoremstyle{definition}
\newtheorem{dfn}[theorem]{Definition}
\newtheorem{example}[theorem]{Example}

\theoremstyle{remark}
\newtheorem{remark}{Remark}[section]
\newtheorem{notat}{Notation}[theorem]
\numberwithin{equation}{section}

\begin{document}

\title[Zeros of polynomials and normal matrices]
{Inverse spectral problem for normal matrices and a
generalization of the Gauss-Lucas theorem}
\author{
S. M. Malamud}

\address{ETH, Z\"urich}

\email{semka\@math.ethz.ch}

\subjclass{Primary: 15A29; Secondary: 30C15, 30C10}

\newcommand{\conv}{\operatorname{conv}}
\newcommand{\gb}{\beta}
\newcommand{\ga}{\alpha}
\newcommand{\Ext}{\operatorname{Ext}}
\newcommand{\CV}{\operatorname{CV}}
\newcommand{\Maj}{\operatorname{Maj}}
\newcommand{\col}{\operatorname{col}}
\newcommand{\tA}{{\widetilde A}}
\newcommand\eps{\varepsilon}
\newcommand\supp{\operatorname{supp}}
\newcommand\diag{\operatorname{diag}}
\newcommand{\cl}{{\mathcal L}}
\newcommand\rank{\operatorname{rank}}
\newcommand\gt{\theta}
\newcommand\gvf{\varphi}
\newcommand\cg{{\mathcal G}}
\newcommand\gve{\varepsilon}
\newcommand\go{\omega}
\newcommand\gO{\Omega}
\newcommand\tr{\operatorname{tr}}
\newcommand\gs{\sigma}
\newcommand\gd{\delta}
\newcommand\gD{\Delta}
\newcommand\D{{\mathbb D}}
\newcommand\Cyc{\operatorname{Cyc}}
\renewcommand\mod{\operatorname{mod}}
\newcommand\Llr{\Longleftrightarrow}
\newcommand\esssup{\operatorname{esssup}}
\newcommand\gS{\Sigma}
\newcommand{\gL}{\Lambda}
\newcommand{\gl}{\lambda}
\newcommand\cd{{\mathcal D}}
\newcommand\C{{\mathbb C}}
\newcommand\R{{\mathbb R}}
\newcommand\N{{\mathbb N}}
\newcommand\wt{\widetilde}
\newcommand{\ca}{{\mathcal A}}
\newcommand{\cb}{{\mathcal B}}
\newcommand{\cs}{{\mathcal S}}
\newcommand{\ch}{{\mathcal H}}
\newcommand{\ck}{{\mathcal K}}
\newcommand{\cm}{{\mathcal M}}
\newcommand{\cx}{{\mathcal X}}

\begin{abstract} We establish an analog of the Cauchy-Poincare
separation theorem  for normal matrices in terms of majorization.
Moreover, we present a solution to the inverse spectral problem
(Borg-type result). Using this result we
essentially generalize and extend the known Gauss--Lucas
theorem about the location of the roots of a complex polynomial and
of its derivative. The last result is applied to prove the
well-known old conjectures of de Bruijn-Springer and Schoenberg.
\end{abstract}

\maketitle

\section{Introduction}

Let $A=A^*$ be a selfadjoint $n\times n$-matrix, $A_{n-1}$ its
principal  $(n-1)\times(n-1)$ submatrix, obtained by deleting the
last row and column. According to the Cauchy--Poincare interlacing
theorem their spectra $\sigma(A)=\{\gl_j\}^n_1$ and
$\sigma(A_{n-1})=\{\mu_j\}^{n-1}_1$ separate  each other, that is
    \begin{equation}\label{0.1}
\gl_1\le\mu_1\le\gl_2\le\mu_2\le \ldots \le\gl_{n-1}
\le\mu_{n-1}\le\gl_n
    \end{equation}
It is known (see \cite{HoJ}, \cite{MarM}) that the converse is
also true, that is for any two sequences $\{\gl_j\}^n_1$ and
$\{\mu_j\}^{n-1}_1$ of real numbers, satisfying \eqref{0.1}, there
exists (nonunique) $n\times n$ selfadjoint  matrix $A$ such that
$\sigma(A)=\{\gl_{j'}\}^n_1$ and
$\sigma(A_{n-1})=\{\mu_j\}^{n-1}_1$. We say that such a  matrix
$A=A^*$ solves the inverse spectral problem for these sequences.

The result of Hochstadt, see \cite{Hoch} (an analog
the well known Borg uniqueness result for Sturm-Liouville
equation) claims that there exists the unique Jacobi (tridiagonal)
selfadjoint matrix $A$ solving the inverse problem.

In the present paper we generalize both the
Cauchy-Poincare separation theorem and the Hochstadt theorem to
the case of normal matrices.


It is obvious that an analogous result can not hold literally
for a normal matrix
$A$,  first because the eigenvalues are not real and second, a
principal submatrix  $A_{n-1}$
is normal only in trivial cases
(see \cite{FanKy} and Lemma \ref{FanKy}).


Our first main result is Theorem \ref{th1.1}.
It provides necessary and sufficient
geometric conditions
for the sequences $\{\gl_j\}^n_1$ and
$\{\mu_j\}^{n-1}_1$
to be the spectrum of a normal matrix $A$
and its submatrix $A_{n-1}$  respectively.


In order to formulate these geometrical conditions we introduce
(in Section 2) several concepts of majorization for sequences of
vectors from $\R^n,$ being  natural generalizations of the
classical ones and coinciding with those for $n=1.$

Note however, that the sufficient part of the above mentioned
conditions can be expressed analytically without
majorization and reads as follows
     \begin{equation}\label{0.2}
c_k:=\frac{\prod_{j=1}^{n-1}(\mu_j-\gl_k)}{\prod_{1\le j\le n,j\not=k}
(\gl_j-\gl_k)}\ge 0,
\qquad  k\in \{1,\ldots, n\}.
   \end{equation}

The second main topic of our paper is concentrated around the known
Gauss-Lucas theorem \cite{PoSe}. According to this theorem the
roots $\{\mu_k\}_1^{n-1}$ of the derivative $p'$ of any complex
polynomial $p(\in \C[z])$ of degree $n$ lie in the convex hull of
the roots $\{\gl_j\}_1^n$  of the polynomial $p$.

At a first view this  topic is rather far from the above one.
Nevertheless the following proposition establishes a "bridge"
between the two main topics (parts) of the paper.
   \begin{prop}\label{01}
Let $p(z)$  be a  polynomial of degree $n$ with zeros
$\{\gl_j\}^n_1$ and  $\{\mu_j\}^{n-1}_1$
the zeros of its derivative.
Then there exists a (nonunique) normal matrix
$A\in M_n({\Bbb C})$ such that $\sigma(A)=\{\gl_j\}^n_1$ and
$\sigma(A_{n-1})=\{\mu_j\}^{n-1}_1$.
     \end{prop}
The proof is immediately implied by \eqref{0.2}
since now $c_k=1,\ k\le n.$

Combining Proposition \ref{01} with Theorem \ref{th1.1} (our
solution to inverse problem) and setting
$\mu_n:=(\sum_1^n\lambda_j)/n,$ we immediately arrive  (see
Proposition \ref{corgau}) at the  following result:
let $\mu:=\{\mu_j\}^n_1,$
$\gl:=\{\gl_j\}^n_1.$  There exists a doubly stochastic
$n\times n$-matrix $S$ such that $\mu= S\gl.$

This result essentially  improves the  Gauss-Lucas theorem.
Proposition \ref{01} allows us to apply linear algebra techniques
to the investigation of the location of
the zeros of a polynomial and of its derivative.
For example, using the exterior algebra techniques we obtain a generalization
of the Gauss-Lucas theorem for the products of roots (see
Proposition \ref{corgau}).

Note, that while the second topic  has attracted a lot of
attention during two last decades (see \cite{Az}, \cite{Cra}, \cite{debr2},
,\cite{Mil}, \cite{Paw}, \cite{Schm}),
our approach seems to be new and perspective. In
particular, in the framework of this approach we get simple (and
short) solutions to the old problems of de Bruijn-Springer
\cite{debr1} and  Schoenberg  \cite{Schoen}.

Let us briefly describe these problems, fixing  the above notations.

In 1948 de Bruijn and Springer \cite{debr1} conjectured that the
following inequality holds for any convex function $f:\C\to\R:$ $$
\frac1{n-1}\sum_{j=1}^{n-1}f(\mu_j)\le
\frac1n\sum_{j=1}^nf(\gl_j). $$ They succeeded in proving this
inequality for a class of convex functions.
We provide a proof of this conjecture by showing that
a bistohastic matrix $S$  in the above mentioned
representation  $\mu= S\gl$ can be choosen in such a way that all entries in the  last row equal $1/n.$
Moreover, we prove (Theorem \ref{debru})
that the
following inequality is valid for any $k\in \{1,\ldots,n-1\}$
      \begin{equation}
\frac1{\binom{n-1}k}\sum_{1\le i_i<\ldots <i_k\le n-1}
f\left(\prod_{j=1}^k(\mu_{i_j}-\ga)\right)
\le \frac{1}{\binom{n}k}\sum_{1\le i_1<\ldots<i_k\le n}
f\left(\prod_{j=1}^k(\gl_{i_j}-\ga)\right).
      \end{equation}

In 1986 Schoenberg \cite{Schoen} (see also \cite{debr2})
conjectured that if $\sum_{j=1}^n\gl_j=0,$ then $$
n\sum_{j=1}^{n-1}|\mu_j|^2\le (n-2)\sum_{j=1}^n|\gl_j|^2 $$ and
the equality holds if and only if all numbers $\gl_j$ lie on the
same line. We establish this inequality in Proposition \ref{01}
(see Theorem \ref{Schoen}).

Let us briefly sketch the contents of the paper.

In section 2 we introduce two new notions of majorization for
sequences of vectors with nonequal numbers of entries and
establish some simple properties of those.
We also study a connection between
different concepts of majorization and show that (on the contrary to
the scalar case), they are not equivalent.
In particular, this provides a negative answer to
the question from  the book
of Marshall and Olkin
\cite{MaOL}, p.433.

In section 3 we establish an analog of
the Cauchy-Poincare separation theorem for normal matrices
(Theorem \ref{th1.1}). As a corollary, we get an analogous
result for "noncommutative" convex combinations of normal
matrices (Corollary \ref{nonc}).
Moreover, we solve an inverse Borg-type problem,
generalizing the result of Hochstadt \cite{Hoch}.

In Section 4
we essentially generalize and extend the known Gauss--Lucas
theorem based on our solution to inverse problem for normal
matrices. Finally, we apply  this result in order to obtain
complete solutions to the de Bruijn-Springer and Schoenberg
conjectures.

The preliminary version of our paper has already been published
as a preprint \cite{Mal4}.
The main results of the paper have been announced without
proofs in \cite{Mal5}.

\begin{center}
{\bf Acknowledgements}
\end{center}
\vspace{3mm}

I would like to express my deep gratitude to M. M. Malamud and F.
V. Petrov for numerous comments, remarks and discussions which
essentially improved the paper. I am also very grateful to J.
Borcea and F. V. Petrov who have informed me about the conjecture
of de Bruijn and Springer \cite{debr1} and the papers \cite{Sher}
and \cite{FiHol}.

\section{Majorization}

{\bf 2.1. Two definitions of majorization. }

We start with several known definitions
(see \cite{Herm}, \cite{MarM}, \cite{MaOL}).
      \begin{notat}
Let $X$ be a subset of $\R^m.$ Denote by $\conv(X)$
the convex hull of $X,$ i.e. the smallest convex set,
containing $X.$

If $X$ is convex, let Ext$X$ denote
the (nonempty) set of its extreme points.

Next, $\CV(Y)$ stands for the set of all convex functions on a
convex set $Y.$

As usual, denote by $C:=A\circ B$  the Schur (element-wise)
product of two $n\times n$-matrices $A=(a_{ij})$ and $B=(b_{ij}):$
$(c_{ij}):=(a_{ij}b_{ij}).$
\end{notat}
     \begin{dfn}
a) A matrix $A\in M_n(\R)$ is called bistochastic
(doubly stochastic) iff all its entries are
nonnegative and the sum of the elements
in each row and each column equals one.

We denote the set of bistochastic  matrices
by $\Omega_n\subset M_n(\R).$

b)  A matrix $A\in M_n(\R)$ is called
unitary-stochastic (orthostochastic)
if there exists  a unitary (orthogonal) matrix
$U\in M_n(\C)(U\in M_n(\R))$
such that  $A=U\circ \bar U.$

The set of all unitary stochastic matrices is
denoted by $\Omega^u_n.$
     \end{dfn}

\begin{remark} The set
of all bistochastic  matrices is convex
and contains all transposition matrices.

The known Theorem of Birkgoff states that the set $\Ext(\Omega_n)$
coincides with the set of all transposition matrices, and thus by
the Krein-Mil'man theorem $\Omega_n$ is the convex hull of all
transposition matrices \cite{Herm}, \cite{MarM}, \cite{MaOL}.

Not that each  unitary stochastic $n\times n$-matrix
is  bistochastic, i.e. $\Omega^u_n\subset \Omega_n.$
But converse is not true: not every bistochastic
matrix is unitary-stochastic \cite{MarM}, \cite{MaOL}.
\end{remark}

         \begin{dfn}\label{defmaj}
Let $x=\{x_k\}_{k=1}^l$ and $y=\{y_k\}_{k=1}^m$
be two sequences of vectors in $\R^n$ and  $l\le m.$
Suppose that the following conditions are fulfilled:
     \begin{equation}\label{2.1}
\begin{split}
&\conv(x_{i_1},\ i_1=\overline{1,l})\subset \conv(y_i:\
i=\overline{1,m}),\\ &\ldots\ldots\ldots\ldots\ldots\\
&\conv(x_{i_1}+\ldots+x_{i_k}:\ 1\le i_1<\ldots<i_k\le l)\subset\\
&\conv(y_{i_1}+\ldots+y_{i_k}:\ 1\le i_1<\ldots<i_k\le m)\\
&\ldots\ldots\ldots\ldots\ldots\\ &x_1+\ldots+x_l\in
\conv(y_{i_1}+\ldots+y_{i_l}:\ 1\le i_1<\ldots<i_l\le m)\\
\end{split}
    \end{equation}
Then we say that the sequence $\{y_k\}$ majorates the sequence
$\{x_k\}$ and write $x\prec y;$
       \end{dfn}
      \begin{remark}
If  $l=m,$ the last condition turns into $$
x_1+\ldots+x_m=y_1+\ldots+y_m. $$
\end{remark}

    \begin{dfn}\label{defbis}
Let $x:=\{x_k\}_{k=1}^l$ and $y:=\{y_k\}_{k=1}^m,$
$l\le m$ be two sets of vectors $x_k,y_k\in \R^n.$

We say  that $x$ is
bistochastically majorated by $y$ and write $x\prec_{ds} y,$
if there exist vectors $x_{l+1},\ldots,x_m\in \R^n$ and a
bistochastic matrix $S\in \Omega_m,$
such that ${\wt x}:=\{x_k\}_1^m=(I_n\otimes S)y.$

If  $S$ can be chosen to be unitary stochastic, we write
$x\prec_{uds}y.$
   \end{dfn}

Next we compare these definitions to the classical
ones.

For this purpose we recall the notion of majorization (\cite{Herm},
\cite{MarM},\cite{MaOL}) for sequences of real numbers
( the case $n=1).$
    \begin{dfn}\label{defbis2}
Let there be given two real sequences $\ga:=\{\ga_k\}_1^m$ and
$\gb:=\{\gb_k\}_1^m.$ Let also $\hat{\ga}$ and $\hat{\gb}$ be
these sequences, reordered to be decreasing. If
      \begin{equation}\label{1.1}
\hat{\gb}_1+\ldots+\hat{\gb}_j\le \hat{\ga}_1+\ldots+\hat{\ga}_j,
\qquad j\in\{1,\ldots,m\}
       \end{equation}
then it is usually written
$\gb\prec\prec\ga.$

If, moreover,
$$
\sum_{k=1}^m\ga_k=\sum_{k=1}^m\gb_k,
$$
then the sequence $\gb$ is said to be majorized by $\ga$
which is denoted by $\gb\prec \ga.$
          \end{dfn}

The following famous theorem due to Weyl, Birkgoff and
Hardy-Littlwood-Polya (see \cite{HLP}, \cite{Herm},
\cite{MarM},\cite{MaOL}), explains the connection between
two different definitions of majorization in the scalar case ($n=1$).

\begin{theorem}\label{major} Let $\ga,\gb\in \R^m$ be to real sequences.
Then the following are equivalent

1) $\gb\prec\ga.$

2) The following inclusion holds true:
        \begin{multline}\label{1.3}
\gb\in \conv(\{A\ga:\ A\text{ is a permutation matrix}\})=\\
\conv(\{\{\ga_{i_1},\ldots,\ga_{i_m}\}:\ \{i_1,\ldots,i_m\}
\text{ is a permutation
of the set }\{1,\ldots,m\}\})
        \end{multline}

3) There exists a bistochastic matrix $S\in M_m(\R),$ such that
$\gb=S\ga.$ It fact, the matrix $S$ can be chosen to be orthostochastic.

4) The inequality
      \begin{equation}\label{1.3'}
\sum_{i=1}^mf(\gb_i)\le \sum_{i=1}^mf(\ga_i)
       \end{equation}
holds for any convex function $f$ on $\R.$

5) The inequality
    \begin{equation}\label{1.3''}
f(\gb_1,\ldots,\gb_m)\le f(\ga_1,\ldots,\ga_m)
     \end{equation}
holds for any convex  function $f$ in $\R^m$ which is symmetric,
that is, invariant under any permutation of the coordinates.
      \end{theorem}
      \begin{remark} \label{rem1}
Note that in the case $n=1$ and $l=m$ definitions \ref{defmaj} and
\ref{defbis} are equivalent to the above Definition \ref{defbis2}
of majorization for real sequences. In fact, we have $$
\gb\prec\ga\Longleftrightarrow \gb=S\ga\Longleftrightarrow
(-\gb)\prec(-\ga). $$ (this fact very easy to check  explicitly).
Now, the convex hull of a set of real numbers is the closed interval
between the minimal to the maximal numbers.
Thus Definition \ref{defmaj} is a natural generalization of the
standard one in $\R^1.$
\end{remark}

It is very easy to see that the following proposition is valid
        \begin{prop}
If $x=\{x_k\}\prec_{ds} y=\{y_k\},$ then $x\prec y.$
      \end{prop}
One can suppose, that a complete  analog of Theorem
\ref{major} is valid, that is the partial orders $\prec$ and $\prec_{ds}$
are equivalent. But the following example shows, that it is not the case.

       \begin{example}\label{exam}
Let $n=2$ and $m=4.$ Set
\begin{equation}
\begin{aligned}
x&=\{x_1,\ldots,x_4\}=\{(12,12),(12,12),(5,3),(3,5)\},\\
y&=\{y_1,\ldots,y_4\}=\{(8,16),(16,8),(0,0),(8,8)\}
\end{aligned}
\end{equation}
It is easy to check by hand, that $x\prec y.$
At the same time it is  easy to see that the vector $x_1=(12,12)$
can be uniquely expressed as a convex combination of $y_k$-s:
$x_1=1/2(y_1+y_2).$
Suppose now, that there exists a bistochastic matrix $S,$ such that $x=Sy.$
Then
$S$ has the form
$$
S=\begin{pmatrix} 1/2&1/2&0&0\\
1/2&1/2&0&0\\
0&0&s_{33}&s_{34}\\
0&0&s_{43}&s_{44}
\end{pmatrix},
$$
which is impossible, since $x_3,x_4$ do not belong to the convex hull of
$y_3,y_4.$
\end{example}

\begin{remark} On the p. 433 of \cite{MaOL}  Marshall and Olkin,
mention the condition 2) of Proposition \ref{ma1} as the weakest
possible notion of majorization. By virtue of Proposition
\ref{ma1}, it is equivalent to $\prec$ for $l=m.$
Furthermore they say that the relation between $\prec$ and $\prec_{ds}$
is not clear. Example \ref{exam} provides a negative answer to
this question.
\end{remark}

We note the following simple

       \begin{prop}
a) Let $\{y_k\}_{k=1}^m$ be such, that $\conv(\{y_k\}_{k=1}^m)$
is affine isomorphic to the standard simplex
$$
\Sigma_{m-1}:= \{(t_1,\ldots,t_m)\in \R^m:\ t_k\ge 0,\ \sum_k t_k=1\}.
$$
Then $x=\{x_k\}\prec y=\{y_k\}$ if and only if
$x\prec_{ds} y.$

b) for $m=3$ the orders $\prec$ and $\prec_{ds}$ are equivalent.
    \end{prop}
      \begin{proof}
a) Under the assumption a) each $x_k$ can be uniquely
represented as a convex combination of $y_k:$
$$
x_k=\sum_{i=1}^ms_{ki}y_i,\qquad k\in \{1,\ldots,m\}
$$
with $\sum_i s_{ki}=1.$ Now, by definition, $x\prec y$ yields
$$
(1/m)\sum_ky_k=(1/m)\sum_k x_k=(1/m)\sum_{i=1}^m\left(\sum_{k=1}^ms_{ki}
\right)y_i.
$$
Defining
$$
\beta_i=1/m\left(\sum_{k=1}^ms_{ki}\right)
$$
we get
$\beta_i=1/m$ for all $i,$ since the expression via extreme points is unique Thus the matrix $S=(s_{ki})_{k,i=1}^m$ is
the required bistochastic matrix.

b) $m=3.$ If the endpoints of $y_1,y_2,y_3$ in $\R^n$
do not lie on the same line, then a) applies.
If they do, then shifting them all by the same vector so,
that the line becomes
passing through the origin we reduce the problem to the case $n=1,$
contained in Theorem \ref{major}.
\end{proof}

\begin{remark} It is clear, that if  $\{y_k\}$ are linearly
independent, then
they satisfy the hypothesis of a).

It is also interesting to note that in the case a) only
the first condition
$x_k\in \conv(\{y_k\}_{k=1}^m)$ and the last one
$\sum x_k=\sum y_k$ are sufficient
for the existence of a bistochastic matrix.
      \end{remark}

{\bf 2.2. The set of extreme points of the set Maj$(y)$.}

Let $n=1$ and $x=\{x_k\}_1^m,\ y=\{y_k\}_1^m\in R^m.$
It is known (see \cite{Mar}), that in this case $(n=1)$
the set of extreme points of the set
Maj$(y):=\{x:\ x\in \R^m,\ x\prec y\}$ is
       \begin{equation}\label{ext}
Ext(\Maj(y))=\{Py:\ P\in \Omega_m,\ P  \text{ is a permutation matrix}\}.
      \end{equation}
The following statement easily follows from the scalar case.

      \begin{prop}
Let $x=\{x_k\}_1^m,\ y=\{y_k\}_1^m$ be two sequences of vectors
in $\R^n$ and Maj$(y):=\{x:\ x\prec y\}.$ Then
$$
Ext(\Maj(y))\supset \{(I_n\otimes P)y:\ P\text{ is a permutation matrix}\}.
$$
      \end{prop}

The following questions naturally arise in this connection:

{\bf Questions:}
1) Find some additional geometric conditions, such that together with
\eqref{2.1} they imply $x\prec_{ds} y.$

2) What are the extreme points of the set $\Maj(y)?$

3) Under which  conditions on the sequence $y=\{y_k\}$ the sets
$\{x:\ x\prec y\}$ and $\{x:\ x\prec_{ds} y\}$ coincide?

Let $CVS(\R^n)$ be the closed in the point-wise convergence topology cone
in $CV(\R^n),$ generated by the set of convex functions
     \begin{equation}\label{om}
\{f(\langle x,y\rangle): f\in \CV(\R),\ y\in \R^n\}.
      \end{equation}

The class  $CVS(\R^n)$ naturally arises in the following
proposition being  a partial generalization of Theorem
\ref{major}.
      \begin{prop}\label{ma1}
Let
$x:=\{x_j\}_1^l,\ y:=\{y_k,\}_1^m$ be systems of
vectors from $\R^n.$
The following conditions are equivalent

1) $x\prec y;$

2) for any vector $h\in \R^n$
$$
(\langle x_1,h\rangle,\ldots,\langle x_l,h\rangle)\prec
(\langle y_1,h\rangle,\ldots,\langle y_m,h\rangle);
$$

3) the inequality
      \begin{equation}\label{2.3''}
\sum_{i=1}^lf(x_i)\le \sum_{i=1}^mf(y_i)
       \end{equation}
holds true for any nonnegative $f\in CVS(\R^n)$ when $l<m$ and
for any $f\in CVS(\R^n)$ if $l=m.$
\end{prop}
       \begin{proof}
1)$\Llr$2). A vector $x\in \R^n$ lies in a convex set
$Y\subset \R^n$ iff its projection to any line
lies in the projection of $Y$
onto the same line.
Thus
$$
x_{i_1}+\ldots+x_{i_k}\in \conv(\{y_{j_1}+\ldots+y_{j_k}\})
$$
iff
$$
\langle x_{i_1},h\rangle+\ldots+\langle x_{i_k},h\rangle\in
\conv(\{\langle y_{j_1},h\rangle +\ldots+\langle y_{j_k},h\rangle \})
$$
for any $h\in \R^n.$ By Remark \ref{rem1} this is equivalent to 2).

2)$\Llr 3).$ By a result of Fisher and Holbrook \cite{FiHol},
$\{\langle x_k,h\rangle\}_1^l\prec \{\langle y_k,h\rangle \}_1^m$
if and only if $$ \sum_{k=1}^l f(\langle x_k,h\rangle)\le
\sum_{k=1}^mf(\langle y_k,h\rangle) $$ for any nonnegative
function $f\in \CV(\R).$ If $l=m,$ we have an equality in
\eqref{2.3''} for any linear function. And any convex function on
the line is a sum of a linear function and a nonnegative convex
function. This immediately yields the required.
\end{proof}

We mention also the following beautiful result, being another partial
generalization of the Hardy-Littlewood-Polya
theorem. It is due to Sherman \cite{Sher}
for $l=m$ and to Fisher and Holbrook \cite{FiHol} for $l<m:$
      \begin{theorem}\cite{Sher}, \cite{FiHol} \label{ThSherman}
Let
$x:=\{x_j\}_1^l$ and  $y:=\{y_k\}_1^m$ be systems of
vectors in $\R^n.$
Then $x\prec_{ds}y$ if and only if
the inequality
      \begin{equation}\label{2.3'}
\sum_{i=1}^lf(x_i)\le \sum_{i=1}^mf(y_i)
       \end{equation}
is valid for any nonnegative $f\in \CV(\C)$ if $l<m$ and for any
$f\in\CV(\C)$ if $l=m.$
    \end{theorem}
The proof is completely different from the scalar case.

Note, that combining Proposition
\ref{ma1}, Theorem \ref{ThSherman} and Example \ref{exam}  we
arrive at the relation $CVS(\R^n)\not =\CV(\R^n).$ Note also that
using our Example \ref{exam} F. V. Petrov has  constructed a
simple explicit counterexample of a function $f\in
\CV(\R^n)\setminus CVS(\R^n).$

To finish the section, we mention the following elegant
result due to F. Petrov.

\begin{prop}\label{cor2.1}
1) $\{x_k\}_{k=1}^l\prec\{y_k\}_{k=1}^m$ if and only if
 the sum of any $s$ of $x_i$-s, $1\le s\le l$ is a linear
combination of $z_j$-s with coefficients between $0$ and $1$ and
the sum of coefficients equal $s;$

2) $\{x_k\}_{k=1}^l\prec\{y_k\}_{k=1}^m$
if and only if there exist vectors $x_{l+1},\ldots,x_m\in \R^n$
such that
$$
\{x_k\}_{k=1}^m\prec\{y_k\}_{k=1}^m.
$$
\end{prop}

\begin{proof} 1) easily follows from Proposition \ref{ma1}, 2).

2) It is easy to see from 1) that if $\{x_k\}_{k=1}^l\prec\{z_k\}_{k=1}^m$,
then $\{x_k\}_{k=1}^{l+1}\prec \{z_k\}_{k=1}^m,$ where
$x_{l+1}={\frac{\sum_1^m z_i-\sum_1^l x_i}{m-l}}.$
\end{proof}

\section{Inverse problem and interlacing theorem for normal
matrices}

In this section we will use the partial orders $\prec$ and
$\prec_{ds}$ for vectors with complex entries. In this case we identify
$\C$ with $\R^2$ and so Definition \ref{defmaj} does not change.

Let $A_{i_1,\ldots,i_k}^{j_1,\ldots,j_k}$ denote the
submatrix of $A$ with the $i_1,\ldots,i_k$ and columns
$j_1,\ldots,j_k.$
We denote for the brevity $A_{n-1}:=A_{1,\ldots,n-1}^{1,\ldots,n-1}.$

 {\bf 3.1. Preliminary solution to the inverse problem by two spectra.}

     \begin{prop}\label{crit}
Let $\{\gl_k\}_1^n$ and $\{\mu_j\}_1^{n-1}$
be two sequences of complex numbers. Then the system of
inequalities:
     \begin{equation}\label{3.1.1}
\frac{\prod_{j=1}^{n-1}(\mu_j-\gl_k)} {\prod_{1\le j\le
n,j\not=k}(\gl_j-\gl_k)}\ge 0, \qquad \  k\in\{1,\ldots,n\}
   \end{equation}
is valid if and only if there exists a normal matrix $A$ with
the spectrum $\gs(A)=\{\gl_1,\ldots,\gl_n\}$ such that the spectrum of
$\tA=A_{n-1}$ is $\gs(\tA)=\{\mu_1,\ldots,\mu_{n-1}\}.$
     \end{prop}
      \begin{proof} $Necessity.$
Let $A$ be a normal matrix with the spectrum
$\gs(A)=\{\gl_1,\ldots,\gl_n\}$ such that
$\gs(\tA)=\{\mu_1,\ldots,\mu_{n-1}\}.$
Let $e=(0,\ldots,0,1).$ Consider the function:
    \begin{equation}\label{3.2}
\gD(\gl):=((A-\gl)^{-1}e,e)=\sum_{k=1}^n\frac{ x_k^2}{\gl_k-\gl}=
\frac{\det(\tA-\gl)}{\det(A-\gl)}
   \end{equation}
where $x_k$ are the coordinates of $e$ in the orthonormal basis
of eigenvectors of $A.$
Clearly, the poles of $\gD(\gl)$ are in the spectrum of $A$ and
the residues in these poles are equal to $x_k^2$ and hence are nonnegative.
But, by \eqref{3.2} these residues equal the numbers \eqref{3.1.1}.

\noindent $Sufficiency.$
Let \eqref{3.1.1} ne fulfilled. Consider the function
     \begin{equation}\label{3.4}
\gD(\gl):=\frac{\prod_{j=1}^{n-1}(\mu_j-\gl)}{\prod_{k=1}^n(\gl_k-\gl)}
     \end{equation}
By \eqref{3.1.1}, the residues of $\gD(\gl)$ in its poles $\gl_k$
are nonnegative  and hence equal $x_k^2$ for some
real numbers $x_k.$ Clearly, we have
    \begin{equation}\label{3.5}
\gD(\gl)=\sum_{k=1}^n\frac{ x_k^2}{\gl_k-\gl}
   \end{equation}
and
$$
\sum_{k=1}^nx_k^2=\lim_{\gl\to \infty} -\gl\gD(\gl)=1.
$$
Therefore, considering the diagonal
matrix $A=\diag\{\gl_1,\ldots,\gl_n\}$ and
writing it down in an orthonormal basis with the last vector
$e_n=(x_1,\ldots,x_n)$ we get the required normal matrix.
       \end{proof}
      \begin{remark}
Note, that the poles of the function  $\gD(\gl)$ (see \eqref{3.2})
are simple (i.e. of multiplicity one), since the matrix $A$ is
(unitary) diagonalizable. Therefore it easily follows from
\eqref{3.2},
as well as from general dimension arguments, that if $A$ has a
$k$-multiple eigenvalue $\gl_0,$ then  $\gl_0$ is an eigenvalue of
$A_{m-1}$ of multiplicity at least $k-1.$ Note, however, that
$A_{n-1}$ is  normal only in very special cases (see Lemma
\ref{FanKy}). Moreover, it may even happen that $A_{n-1}$ is not
of simple structure, that is it may be nondiagonalizable.
      \end{remark}

\begin{cor}\label{crit2} Let two systems of complex numbers
$\{\mu_j\}_1^{n-1}$ and $\{\gl_j\}_1^n$ satisfy
      \begin{equation}
\Delta(\gl):=\frac{\prod_{j=1}^{n-1}(\mu_j-\gl)}
{\prod_{j=1}^{n}(\gl_j-\gl)}=\sum_{j=1}^n\frac{
|x_j|^2}{\gl_j-\gl}
     \end{equation}
with some complex numbers $\{x_j\}_1^n.$
Then for any unitary matrix
$U=(u_{ij})_{i,j=1}^n$ with the last row
$(u_{n1},\ldots, u_{nn})=(x_1,\ldots, x_n),$
the matrix $A:=U\diag\{\gl_j\}_1^nU^*$ satisfies the
hypothesis of Proposition \ref{crit}.
     \end{cor}

{\bf 3.2.  Quasi-Jacobi normal matrices and an analog
of the  Hochstadt theorem.}

It is known and easy to see that any selfadjoint matrix
(bounded operator) is unitary equivalent to a selfadjoint
tridiagonal (Jacobi) matrix.

Here we  find an analog of such a  form for a normal
matrix and apply it in order to obtain an  analog of the
Hohstadt result \cite{Hoch} on  the unique recovery of a Jacobi
matrix from two spectra.


\begin{prop} \label{jac}
Every normal $m\times m$ matrix $A$ is unitary equivalent to a direct
sum of normal matrices $A_i,\ i=1,\ldots,k$ satisfying
$(A_i)_{j,k}=0$ for $k\ge j+2$ and $(A_i)_{j,j+1}\not=0.$
Moreover, $A$ has simple spectrum iff it is unitary equivalent to
only one such matrix.
\end{prop}

\begin{proof}
The proof is very simple and standard. It is clear that it
suffices to consider only the case of simple spectrum.

In this case, taking any cyclic vector $x,$ we get that
$\{A^jx\}_0^{m-1}$ forms  a basis in $\C^m.$ After the
Gram-Schmidt procedure we arrive at the required basis.
\end{proof}

\begin{dfn}
A matrix $A$ is called quasi-Jacobi if it satisfies the hypothesis
of Proposition \ref{jac}.
\end{dfn}

Thus, quasi-Jacobi form is in a sense a normal form for a normal
matrix. Now we can complement Proposition \ref{crit} with a
uniqueness result.

      \begin{theorem}
For any two systems of complex numbers
$\{\gl_j\}_1^n$ and $\{\mu_j\}_1^{n-1},$
satisfying \eqref{3.1.1} there exists a unique normal quasi-Jacobi
matrix $A$ such that
$\sigma(A)=\{\gl_j\}_1^n$ and $\sigma(A_{n-1})=\{\mu_j\}_1^{n-1}.$
 \end{theorem}

\begin{proof} Writing the function $\Delta(\gl)$ from
\eqref{3.4} in the form
$$
 \Delta(\gl)=\int\frac{d\mu(z)}{\gl-z}
$$
with $d\mu=\sum_{k=1}^n x_k^2\delta_{\gl_k}$ we  can introduce
the orthogonal polynomials with respect to measure just like in
the Jacobi case (see, e.g. \cite{Ach} and Gesztesy and Simon
\cite{GeSim}) and then, following the same lines as in
\cite{GeSim}, we get the result.
\end{proof}

\begin{remark}
The function $\Delta(\gl)$ is an analog of the Weyl M-function in
this case (see \cite{GeSim}, \cite{Mal1}). Note that our proof of
Proposition \ref{crit} is similar to that proposed in
\cite{GeSim}, \cite{Mal1}.
     \end{remark}

{\bf 3.3. The set of all possible diagonals in the "unitary" orbit
of a  normal matrix.}

The criterion \eqref{3.1.1} of Proposition \ref{crit}
is trivial and provides no information on the geometry
of the sequences $\{\gl_k\}_1^n$ and $\{\mu_j\}_1^{n-1}.$
It is even unclear how far can $\mu_j$ lie from $\gl_k.$
Therefore it would
be desirable to have a more "geometric" answer, being an
analog of the Poincare Theorem.

In this subsection we start with an arbitrary normal matrix  $A$
and give a (rather trivial) description of the set of diagonals of
its "unitary" orbit $\{ UAU^*:\ U\in M_n(\C),\ U^*U=I.\}$ In the
next section we apply this result to complete solution to the
inverse spectral problem for a normal matrix.

         \begin{prop}\label{prop2.1}
Let $A\in M_n(\C)$ be a
normal matrix with the spectrum
$(\gl_1,\ldots,\gl_n).$ Then

a) $(a_{11},\ldots,a_{nn})\prec_{uds} (\gl_1,\ldots,\gl_n);$

b) there exists an orthonormal basis $\{e_i\}_1^n$ such that
$(Ae_i,e_i)=\ga_i$ iff there exists a unitary stochastic  matrix
$O$ such that $\col(\ga_1,\ldots,\ga_n)=O\col(\gl_1,\ldots\gl_n).$

c) if $A$ is selfadjoint, then the set of all possible
diagonals (in all orthonormal bases) is convex;

d) the set of all possible diagonals (in all orthonormal bases)
of a fixed normal matrix $A$ is not necessarily convex.
\end{prop}

\begin{proof} The validity of a) and b) is obvious.
It is also clear that all permutations of the set $(\gl_1,\ldots,\gl_n)$
are realized by diagonals. Thus, if the set of diagonals were convex,
it would
contain all vectors of the form $S\col(\gl_1,\ldots,\gl_n)$ with
a bistochastic $S$.

c) is due to Horn \cite{HoJ}.

d) take $$ S=\frac12
\begin{pmatrix}
1&1&0\\
1&0&1\\
0&1&1
\end{pmatrix}
$$
It is known (see \cite{MarM}) that $S$ is not unitary-stochastic.
Set
$(\gl_1,\gl_2,\gl_3)=(1,i,0)$ and
$$
S(\gl_1,\gl_2,\gl_3)=(\ga_1,\ga_2,\ga_3)^t=(1+i,1,i)^t.
$$
If there exists a unitary-stochastic matrix $O,$ such that
$$
O(\gl_1,\gl_2,\gl_3)^t=(\ga_1,\ga_2,\ga_3)^t,
$$
then it immediately yields  $O=S.$
      \end{proof}

 \begin{remark}
The result of c) is a due to Horn \cite{MarM}, \cite{MaOL},
\cite{Herm}, \cite{HoJ}. All other statements are folklore.
    \end{remark}

{\bf 3.4. Analog of the Cauchy-Poincare interlacing theorem and
a solution to the inverse spectral problem for normal
matrices.}

Now we are ready to state the main result of the section.

Define for any vector $\{\gl_j\}_1^m\in \C^m$ the vector
\begin{equation}\label{ext1}
C_k(\{\gl_j\}_1^m):=\{\gl_{i_1}\cdots\gl_{i_k}\}_{1\le
i_1<\ldots<i_k\le m}\in C^{\binom{m}{k}}.
\end{equation}
      \begin{theorem} \label{th1.1}
Let $\{\gl_1,\ldots,\gl_n\}$ and $\{\mu_1,\ldots,\mu_{n-1}\}$ be
two systems of complex numbers. Then for the existence of a normal
matrix $A$ such that $\gs(A)=\{\gl_1,\ldots,\gl_n\}$ and
$\gs(A_{1,n-1})=\{\mu_1,\ldots,\mu_{n-1}\}$ it is necessary that
the condition
    \begin{equation}\label{3.6}
C_k(\{\mu_j-\ga\}_1^{n-1})\prec_{uds} C_k(\{\gl_j-\ga\}_1^n)
      \end{equation}
be fulfilled for any complex number $\ga\in \C$ and any
$k\in\{1,\ldots, n-1\}$ and sufficient that
it be fulfilled for $k=n-1$ and all
$\ga\in \{\gl_1,\ldots,\gl_k\}.$
      \end{theorem}
       \begin{proof}
a) {\bf Sufficiency.} Let $\ga=\gl_k.$ \eqref{3.6} for $k=n-1$
reads
\begin{multline} \prod_{i=1}^{n-1}(\mu_i-\gl_k)\in \conv
\left(\prod_{1\le i\le n, i\not=l}\{\gl_i-\gl_k\}:\
l=\overline{1,n}\right)\\
= \conv \left(0;\ \prod_{1\le i\le
n,i\not=k}\{\gl_i-\gl_k\}\right)
\end{multline}
Hence \eqref{3.1.1} is valid for all $k.$ Proposition \ref{crit}
yields the required.

b) {\bf Necessity.} Let $A$ be a normal matrix with the spectrum
$\gs(A)=\{\gl_1,\ldots,\gl_n\}$ and $\gs(A_{n-1})
=\{\mu_1,\ldots,\mu_{n-1}\}.$ Let us prove \eqref{3.6}. By the
Shur Theorem \cite{MarM}, \cite{HoJ}
there exists a unitary matrix $V_1\in M_{n-1}(\C)$
such that the matrix $V_1^*A_{1,n-1}V_1$ is upper triangular.
Therefore, considering the matrix $U_1:=V_1\oplus 1\in M_n(\C)$ we
get the normal matrix $B:=U_1^*AU_1$ with the same spectrum as
$A$, $\sigma(B)=\sigma(A),$ but the $\mu_j$-s are on the diagonal.
Therefore we can take $B$ instead of $A.$
Proposition \ref{prop2.1} implies $\{\mu_j\}_1^{n-1}\prec_{uds}
\{\gl_j\}_1^n.$

Take an arbitrary $\ga\in\C$. The matrix $B-\ga I$ is also normal.
Let us consider its exterior power $C_k(B-\ga I):=\wedge^k(B-\ga
I)$ acting on the space$
\wedge^kH:=\underbrace{H\wedge\cdots\wedge H}_{k\text{ times}}.$
Then this matrix is also normal with the spectrum $$ \gs(C_k(B-\ga
I))= \{(\gl_{i_1}-\ga)\cdots (\gl_{i_k}-\ga)\}_{1\le
i_1<\ldots<i_k\le n}. $$ The diagonal elements of $C_k(B-\ga I)$
are the $k\times k$ principal minors of $B-\ga I.$ Since
$B_{n-1}-\ga I_{n-1}$ is upper triangular, then the $k\times k$
principal minors of $B_{n-1}-\ga I$ equal
$$(\mu_{i_1}-\ga)\ldots(\mu_{i_k}-\ga), 1\le i_1<\ldots<i_k\le
n-1.$$ Thus Proposition \ref{prop2.1}, a) implies the required.
      \end{proof}

\begin{remark}
It is easy to see from the proof that the unitary stochastic
matrices in \eqref{3.6} are independent of $\ga.$
\end{remark}

      \begin{cor} There exists
a normal matrix $A$ such that $\gs(A)=\{\gl_1,\ldots,\gl_n\}$ and
$\gs(A_{n-1})=\{\mu_1,\ldots,\mu_{n-1}\}$ if and only if
$$
\prod_{k=1}^{n-1}(\gl_j-\mu_k) \in \conv(0;\ p'(\gl_j))\  \text{
for all }\  j\in \{1,\ldots,n\}
$$
where $p(\gl)=\prod_{k=1}^n(\gl-\gl_k)$
      \end{cor}
       \begin{example}
In the case $n=3$ the orders $\prec$ and $\prec_{ds}$ are
equivalent. Therefore in this case conditions \eqref{3.6} take a
specially simple form:
\begin{align*}
&\mu_1,\mu_2\in\conv\{\gl_j\}_1^3,\ \mu_1+\mu_2\in \conv
\{\gl_1+\gl_2,\gl_2+\gl_3,\gl_1+\gl_3\},\\
&\mu_1\mu_2-\ga(\mu_1+\mu_2)\in
\conv\{\gl_k\gl_p-\ga(\gl_k+\gl_p)\}_ {1\le k<p\le 3},\ \ga\in \C.
\end{align*}
       \end{example}

An immediate consequence is:
       \begin{cor} \label{pos}
Let $A$ be a normal matrix with the spectrum $\{\gl_1,\ldots,\gl_n\}$
and $P$ an $m(\le n)-$dimensional orthoprojection in $\C^n.$ Let
also $B:=PA\lceil P\C^n$ and $\gs(B):=(\mu_1,\ldots,\mu_m).$ Then
    \begin{equation}\label{3.6.1}
C_k(\{\mu_j-\ga\}_1^{m})\prec_{uds} C_k(\{\gl_j-\ga\}_1^n)
      \end{equation}
for all $k,\ 1\le k\le m.$
        \end{cor}
         \begin{cor}\label{nonc}
Let $A_i\in M_{n_i}(\C),\ i=1,\ldots,p$ be a $p$-tuple of normal
matrices, and $\gs(A_i)=\{\gl^i_k\}_{k=1}^{n_i}$ their spectra.
Let $S_i$ be $m\times n_i$ matrices, $i=\overline{1,p}$ such
that
   \begin{equation}\label{4.1.1}
\sum_{i=1}^p S_i^*S_i=I_m
    \end{equation}
is the identity operator in $\C^m.$ Consider
     \begin{equation}\label{4.1.2}
B=\sum_{i=1}^p S_i^*A_i S_i
     \end{equation}
and let $\gs(B)=(\mu_1,\ldots,\mu_m).$ Set $n=n_1+\ldots+n_p$ and
    \begin{equation}\label{4.1}
(\gl_1,\ldots,\gl_n)=(\gl^1_1,\ldots,\gl^1_{n_1},\gl^2_1,\ldots,
\gl^2_{n_2},\ldots,\gl^p_{n_p}).
      \end{equation}
Then the systems of numbers $(\gl_1,\ldots,\gl_n)$ and
$(\mu_1,\ldots,\mu_m)$ satisfy conditions \eqref{3.6.1}.
\end{cor}

\begin{proof} Consider the normal matrix $A:=\oplus_{i=1}^l A_i.$
Condition \eqref{4.1.1} means that the operator
$$
V:=\begin{pmatrix}S_1\\ \vdots \\ S_n
\end{pmatrix}:\ \R^m\to \R^n$$
is an isometry and hence
the operator
$$
B=V^*AV
$$
is unitary equivalent to $PA\lceil P\C^n$ where $P$ is the
projector onto the image of $V.$ Corollary  \ref{pos} completes the
proof.
\end{proof}

\begin{remark} There is another way to implement the trick of going from
Corollary \ref{pos} to Corollary \ref{nonc}. See \cite{Mal2},
\cite{Mal3}.
        \end{remark}

\section{Location of roots of a polynomial and of its
derivative}

{\bf 4.1. Generalization of the Gauss-Lucas Theorem.}

Recall  the known Gauss-Lucas theorem
      \begin{theorem}
The roots of the  derivative $p'$ of a polynomial $p\in C[z]$ lie
in the convex hull of the roots of $p$.
      \end{theorem}
Numerous papers are devoted to different generalizations and
improvements of this results (see e.g. \cite{Dim}, \cite{Az},
\cite{Schm}).

In what follows we denote by $p\in C[z]$ a  polynomial
of degree $n$ with  complex coefficients. Let also $\{\lambda_j\}^n_1$ be the
roots of $p$ and $\{\mu_j\}^{n-1}_1$ the roots of its derivative
$p'$. We set additionally $\mu_n:=(\sum_1^n\lambda_j)/n.$ Then the
Gauss-Lucas theorem reads as follows: there exists  a stochastic
(by rows) matrix $S'$ such that
$\mu=S'\gl,$ where
$\mu:=\{\mu_j\}_1^n$ and $\gl:=\{\gl_j\}_1^n.$

The following result, being a corollary of Theorem  \ref{th1.1},
improves the Gauss-Lucas theorem.
       \begin{prop}\label{corgau}
Let $p(\in \C[z])$ be a degree $n$ complex polynomial with roots
 $\{\gl_j\}_1^n$ and $\{\mu_k\}_1^{n-1}$ the
roots of its derivative $p'$. Then the the sequences
 $\{\gl_j\}_1^n$ and $\{\mu_k\}_1^{n-1}$  satisfy \eqref{3.6}
for $k\in \{1,\ldots,n-1\}$. In particular, there exists a matrix
$S\in \Omega_n$ such that $\mu=S\gl,$ that is, the vector $\mu$ is
bistochastically majorized by  the vector  $\gl.$
         \end{prop}
      \begin{proof} It follows from
Corollary \ref{crit2} and the obvious identity
      \begin{equation}\label{111}
p'(z)/p(z)=\sum1/(z-\lambda_j)
     \end{equation}
that there exists a normal matrix $A$ such that
$\sigma(A)=\{\gl_j\}_1^n$ and $\sigma(A_{n-1})=\{\mu_j\}_1^{n-1}.$
One completes the proof by applying Theorem \ref{th1.1}.
       \end{proof}
Note that this corollary does not give a complete information on the
location. In particular, it does not anyhow explain the
identities:

\begin{equation}\label{prodeq}\frac{1}{C_{n-1}^k}\sum_{1\le i_1<\ldots<i_k\le
n-1}\prod_{j=1}^k(\mu_{i_j}-\ga)=\frac{1}{C_{n}^k}\sum_{1\le i_1<\ldots<i_k\le
n}\prod_{j=1}^k(\gl_{i_j}-\ga)
\end{equation}

which mean that even the products of the roots are equally
distributed.

It turns out that in this case we can obtain a more complete
information on the matrices $S\in \Omega^u_n,$ realizing
majorization.
      \begin{theorem} \label{10}
Let $k\in\{1,\ldots,n-1\},$ $a:=\binom{n-1}k,\ b:=\binom{n}k.$
Then there exists a matrix $S_k=(s_{ijk})\in \Omega^u_b$ such that
$\sum_{i=a+1}^{b}s_{ijk}=k/n$ for all $j\in\{1,\ldots, b\}$ and
$C_k(\{\mu_j-\ga\}_1^{n-1})=PS_kC_k(\{\gl_j-\ga\}_1^{n})$ (see
\eqref{ext1}), where $P:\C^{b}\to\C^{a}$ is the natural
orthoprojection.

In particular, if $k=1$  then there exists $S_1\in \Omega^u_n$
such that $s_{nj1}=1/n,\ j\in \{1,\ldots,n\}$ and
$\{\mu_j\}_1^{n}=S_1\{\lambda_j\}_1^n.$
       \end{theorem}

\begin{proof}
Set $D:=\diag(\gl_1,\ldots,\gl_n).$  By Corollary \ref{crit2} and
\eqref{111} we have, that for any unitary matrix $V$ with the last
row consisting of $1/\sqrt{n},$
the normal matrix $A:=VDV^*$ solves the inverse problem
for a pair of sequences $\{\lambda_j\}_1^n$ and
$\{\mu_j\}_1^{n-1},$ that is  $\sigma(A)=\{\lambda_j\}_1^n$
and $\sigma(A_{n-1})=\{\mu_j\}_1^{n-1}.$
Let, further $U_1$ be the same unitary matrix as in the proof
of Theorem \ref{th1.1}. Then $U:=U_1V(\in M_n(\C))$
is a unitary matrix
with the last row consisting of $1/\sqrt{n}$.
Moreover, it follows from the proof of
Theorem \ref{th1.1} that
$\{\mu_j\}_1^{n}=S_1\{\lambda_j\}_1^n,$
where $S_1:=U\circ \bar U(\in \Omega^u_n)$
is the unitary stochastic matrix
with the last row consisting of $1/n.$

Therefore, passing to the exterior powers as in the proof of
Theorem \ref{th1.1}, we get that the unitary stochastic matrix $$
S_k:=C_k(U)\circ C_k(\bar U), \qquad k\in \{1,\ldots,n\} $$
realizes the unitary stochastic majorization of the systems of
numbers $C_k(\{\mu_j\}_1^{n-1})$ and $C_k(\{\gl_j\}_1^n)'.$ We
prove that $S_k$ has additional properties
          \begin{equation}\label{ident}
\sum_{1\le i_1<\ldots< i_{k-1}\le
n-1}\left|U_{i_1,\ldots,i_{k-1},n}^{1,\ldots,k}\right|^2=\frac{n-k}n.
          \end{equation}
Because of the symmetry, the same identity is certainly valid
for all other choices of $k$ columns.

Expanding each minor with respect to the last row we get
\begin{equation}
U_{i_1,\ldots,i_{k-1},n}^{1,\ldots,k}=\frac1{\sqrt{n}}\sum_{l=1}^k(-1)^{l+1}
U_{i_1,\ldots,i_{k-1}}^{1,\ldots,\hat{l},\ldots,k}
\end{equation}
and thus
     \begin{equation}\label{lapl}
|U_{i_1,\ldots,i_{k-1},n}^{1,\ldots,k}|^2=\frac1n
\sum_{p,q=1}^kU_{i_1,\ldots,i_{k-1}}^{1,\ldots,\hat{p},k}
{\bar{U}}_{i_1,\ldots,i_{k-1}}^{1,\ldots,\hat{q},\ldots,k}
      \end{equation}
Consider now the matrix $V:=U_{1,\ldots,n-1}^{1,\ldots,k}.$ Then,
since $U$ is unitary, we get
$$
T:=V^*V=\{t_{ij}\}_{i,j=1}^k
$$
with
$$
t_{ij}=
\begin{cases}
\frac{n-1}n, & i=j\\
-\frac1n, & i\not=j
\end{cases}.
$$ Therefore we get $$ C_{k-1}(V)^*\cdot C_{k-1}(V)=C_{k-1}(T). $$
Because of \eqref{lapl}, one easily sees that the required sum in
\eqref{ident} is a linear combination of the elements of the
matrix $C_{k-1}(V)^*\cdot C_{k-1}(V)$. Thus this sum depends only
on the matrix $T,$ is independent of the choice of $U$ and is the
same for all other choices of $k$ columns.

But by \eqref{prodeq}, these sums can be nothing but $(n-k)/n.$
The proof is complete.
\end{proof}

      \begin{cor} Let $U$ be an $n\times n$ unitary matrix with
$u_{nj}=1/{\sqrt{n}}.$ Let  $C_k(U)$ be its $k$-th exterior power.
Then for $k\le n$ $S_k:=C_k(U)\circ C_k(\bar U)$ is a unitary
stochastic matrix, satisfying $$
\sum_{i=1}^{\binom{n-1}{k}}s_{ij}=\frac{\binom{n-1}{k}}{\binom{n}{k}}=
\frac{n-k}{n} $$ for all $j\in\{1,\ldots, \binom{n}{k}\}.$
      \end{cor}

{\bf 4.2. Conjecture of de Bruijn-Springer and its
generalization.}

In 1947 de Bruijn and Springer \cite{debr1} conjectured that the
inequality
      \begin{equation}\label{bruin1}
\frac1{n-1}\sum_{j=1}^{n-1}f(\mu_j)\le \frac1n\sum_{j=1}^nf(\gl_j)
      \end{equation}
holds for any convex function $f:\C\to\R.$

In order to prove  this conjecture as well as its generalization
we need the following simple lemma.
      \begin{lemma}\label{conlemma}
\label{triv} Let $\{x_j\}_1^{k},\ \{y_j\}_1^n$ be two sequences of
vectors from $\R^m.$ Suppose that there exists a matrix
$S=(s_{ij})\in \R^{k\times n}$ with non-negative entries and such
that  $\{x_j\}_1^{k}=(S\otimes I_m)\{y_j\}_1^n$ and
     \begin{equation}\label{5.1A}
\sum_{j=1}^ns_{ij}=1, \ i\in \{1,..,k\}\quad  \text{ and }\quad
\sum_{i=1}^ks_{ij}=\frac{k}{n}, \  j\in\{1,...,n\}.
       \end{equation}
Then the inequality
        \begin{equation}\label{5.1B}
\frac1k\sum_{j=1}^k f(x_j)\le \frac1n\sum_{j=1}^n f(y_j)
      \end{equation}
holds true for any function $f\in CV(\R^m).$
       \end{lemma}
      \begin{proof}
One obtains the proof by combining the Jensen inequality with
relations \eqref{5.1A}.
        \end{proof}
The following result contains in particular
a positive solution to the conjecture of
de Bruijn and Springer \cite{debr1}.
           \begin{theorem} \label{debru}
The following inequality holds
true for any convex function $f:\C\to
\R$ and any $k,\ 1\le k \le n:$
        \begin{equation}\label{bruin}
\frac1{\binom{n-1}k}\sum_{1\le i_i<\ldots <i_k\le
n-1}f\left(\prod_{j=1}^k(\mu_{i_j}-\ga)\right)\le \frac{1}{\binom{n}k}
\sum_{1\le i_1<\ldots<i_k\le
n}f\left(\prod_{j=1}^k(\gl_{i_j}-\ga)\right).
         \end{equation}
             \end{theorem}
       \begin{proof}
The inequality immediately follows by combining
Theorem \ref{10} with Lemma \ref{conlemma}.
      \end{proof}
     \begin{remark}
In the case $k=1$ inequality \eqref{bruin} coincides with
inequality \eqref{bruin1}, that is with the
de Bruijn--Springer conjecture \cite{debr1}.
      \end{remark}
        \begin{remark}
Acording to the result of Sherman \cite{Sher}, the existence of a
$k\times n$ matrix $S,$ satisfying the hypothesis of Lemma
\ref{conlemma} is actually equivalent to the validity of
inequality \eqref{5.1B} for each function $f\in CV(\R^m).$
        \end{remark}

{\bf 4.3. The Schoenberg conjecture.}

Now we are ready to prove
the famous Schoenberg conjecture \cite{Schoen},\cite{debr2}.

We will need two Lemmas. The first one is known
\cite{MarM}, but we present it with a proof for the reader's
convenience.

\begin{lemma}\label{HS} Any matrix $A=(a_{ij})_{i,j=1}^n
\in M_n(\C)$ with spectrum
$\sigma(A)=\{\gl_j\}_1^n$ satisfies the inequality
       \begin{equation}\label{00}
\sum_{j=1}^n|\gl_j|^2\le \|A\|_2^2=\sum_{i,j=1}^n|a_{ij}|^2
       \end{equation}
and the equality holds if and only if $A$ is normal.
\end{lemma}

\begin{proof} The inequality \eqref{00} is known.
It is clear that the equality holds true for a normal matrix.

Conversely, let $A$ satisfy the equality. By the Schur theorem $A$ is
unitary equivalent to an upper triangular matrix with $\gl_j$-s
on the diagonal. Since $\|A\|_2$
is unitary invariant, this matrix will be diagonal, that is
$A$ is normal.
       \end{proof}
       \begin{lemma}\label{HSA}
Let $\eps=e^{2\pi i/n}$ and $U=n^{-1/2}(\eps^{k(j-1)})_{k,j=1}^n.$
Let also $\sum_{j=1}^n\gl_j=0.$ Define $$
r(z):=\sum_{j=1}^n\gl_jz^{j-1} $$ and $$
A:=U\diag(\gl_j)_{j=1}^nU^*=\frac1n(r(\eps^{k-j}))_{k,j=1}^n=:(a_{ij})_{i,j=1}^n.
$$ Then $A$ is a normal matrix with spectrum
$\sigma(A)=\{\gl_j\}_1^n$ and $\sigma(A_{n-1})=\{\mu_j\}_1^{n-1}.$
Moreover the following identity holds true $$
n\|A_{n-1}\|_2^2=(n-2)\|A\|_2^2. $$
\end{lemma}
       \begin{proof}
The first statement follows from Corollary \ref{crit2},
since the last row of $U$ consists of $\frac1{\sqrt{n}}.$
It remains to
prove the last identity.

It is easy to see that $a_{jj}=0$ and
$a_{nj}=a_{n-k,j-k}$ for $k<j$
and $a_{nj}=a_{k-j,k}$ for $k>j.$
Threfore the required identity takes the form
$$
\sum_{j=1}^n|\gl_j|^2=n\sum_{j=1}^{n-1}|a_{nj}|^2.
$$
But
$$
n\sum_{j=1}^{n-1}|a_{nj}|^2=\frac{(n-1)}n\sum_{j=1}^n|\gl_j|^2+
\frac1n2\sum_{i<j}\Re\left(\gl_i\bar{\gl_j}c_{ij}\right)
$$
where
$$ c_{ij}=\sum_{k=1}^{n-1}\eps^{(i-j)(n-k)}=-1$$ for all $1\le
i\not=j\le n.$ But we have
$$
\left|\sum_{j=1}^n\gl_j\right|^2=0\Llr2\sum_{i<j}\Re(\gl_i\bar{\gl_j})
=-\sum_{j=1}^n|\gl_j|^2.
$$
This completes the proof.
       \end{proof}
The next lemma is due to Fan Ky and Pall \cite{FanKy}.
      \begin{lemma}\label{FanKy}
Let $A$ be a normal matrix such that its submatrix $A_{n-1}$ is
also normal and $A\not=A_{n-1}\oplus a_{nn}.$ Then all the
eigenvalues of $A$ lie on the same line.
       \end{lemma}
The following result has been conjectured by Schoenberg
\cite{Schoen} (see also \cite{debr2}).
          \begin{theorem} \label{Schoen}
Let $\sum_{j=1}^n\gl_j=0.$ Then
$$
n\sum_{j=1}^{n-1}|\mu_j|^2\le (n-2)\sum_{j=1}^n|\gl_j|^2
$$
and the equality holds if and only if all the numbers $\gl_j$
lie on the same line.
      \end{theorem}
      \begin{proof}
Combining Lemma \ref{HS}  with Lemma  \ref{HSA}, we get
 \begin{equation}\label{last}
n\sum_{j=1}^{n-1}|\mu_j|^2\le n\|A_{n-1}\|_2^2=
(n-2)\|A\|_2^2=(n-2)\sum_{j=1}^n|\gl_j|^2.
    \end{equation}
Moreover, by Lemma \ref{HS}  identity \eqref{last}
holds if and only if  $A_{n-1}$ is normal.
On the other hand, by Lemma \ref{FanKy} this is
possible if and only if all $\gl_j$-s lie on the same line.
\end{proof}

{\bf 4.4. The Mason-Shapiro Polynomials.} In \cite{MS} Gisli
Masson and Boris Shapiro initiated study of a class of
differential operators $T_Q$ defined as follows: let $Q$ be a
degree $k$ monic poplynomial. Then $T_Q$ is defined via $$ T_Q:\
f\to (Qf)^{(k)}. $$ They have shown that for each $m$
there exists a unique polynomial eigenfunction $p_m$ of $T_Q$ of
degree $m.$ Moreover,
     \begin{equation}\label{Shap1}
T_Qp_m=\gl_{m,k}p_m
     \end{equation}
and $\gl_{m,k}$ depends only on $k,m,$ namely
$$\gl_{m,k}=(m+1)(m+2)\cdots(m+k).$$

One of the results of \cite{MS} is  the following interesting analog of the
Gauss-Lucas Theorem:

\begin{theorem}\label{Shap} The zeros of $p_m$ are contained in
the convex hull of the set of zeros of $Q$ for each $m.$
\end{theorem}

The authors have also made a number of beauteful conjectures
about the asymptotic
distribution of the zeros of $p_m$-s, recently proved in
\cite{Berg}.

We strengthen Theorem \ref{Shap} in the following way:

\begin{theorem} Let $\{z_j\}_{j=1}^k$ be the zeros of $Q$ and
$\{w_j\}_1^m$ the zeros of $p_m.$ Then

(i) there exists an $m\times k$
matrix $S$ with
$\sum_{j=1}^ks_{ij}=1$ for all $i\in \{1,\ldots,m\}$ and
$\sum_{i=1}^m s_{ij}=m/k$ for all
$j\in \{1,\ldots,k\}$ and such that
   \begin{equation} \label{lastmatrix}
\{w_j\}_1^m=S\{z_j\}_1^k;
    \end{equation}

(ii) for any convex
function $f:\C\to \R:$
    \begin{equation} \label{HarShap}
\frac{\sum_{j=1}^mf(w_j)}m\le
\frac{\sum_{j=1}^kf(z_j)}k.
      \end{equation}
\end{theorem}
     \begin{proof}
Applying  Theorem \ref{10} $k$ times to the polynomial $Qp_m$
and its consequentive derivatives we arrive at the representation
   \begin{equation}\label{shaprepres}
\{w_j\}_1^m=\widetilde{S}(\col(\{w_j\}_1^m,\{z_j\}_1^k))
     \end{equation}
with $m\times (m+k)$ matrix $\widetilde(S)=(s'_{ij})$
being a product of the $k$ corresponding matrices and
satisfying
$$
\sum_{j=1}^ks'_{ij}=1,\  i\in \{1,\ldots,m\}\  \text{and} \
\sum_{i=1}^m s'_{ij}=\frac{m}{m+k}, \ j\in \{1,\ldots,m+k\}.
$$
Now applying Lemma \ref{conlemma}
we arive at the inequality
$$
\frac1m\sum_{j=1}^mf(w_j)\le
\frac1{k+m}\left(\sum_{j=1}^kf(z_j)+\sum_{j=1}^mf(w_j)\right)
$$
which yields \eqref{HarShap}.

It is
not diffifult to construct the matrix $S$ satisfying \eqref{lastmatrix}
and  the  other required properties
 by resolving the last
identity for $w_j$-s in \eqref{shaprepres}.
But its existence
is immediately implied  by \eqref{HarShap}
due to the result of Sherman \cite{Sher} (see Theorem \ref{ThSherman}).
        \end{proof}

\begin{remark} In the case of real numbers $w_j$ and $z_j$
inequality \eqref{HarShap} has been mentioned (without proof)
by  Harold Shapiro \cite{Shap}.
\end{remark}

{\bf Final remarks.} 1) Let $A\in M_n(\C)$ be a normal matrix, $e\in\C^n$
a vector and $P$ the orthoprojection onto the orthogonal complement of $e.$
It is not difficult to construct examples of a nondiagonalizable
$A_{e}:=PA\lceil P\C^n.$
It would be interesting to investigate the Jordan structure and
other  similarity (or unitary) invariants of the operator $A_{e}.$

2) We do not know whether the relations \eqref{3.6.1} are also sufficient.

3) The most famous unsolved conjectures connected with the
Gauss-Lucas Theorem are the conjectures of Sendov and Smale (see
\cite{Schm} for a survey on this topic).

4) the relative location of the zeros of $p(z),\ p'(z)$ and $p''(z)$ may be very
nontrivial even in the case of real roots \cite{ShaSha}.

5)  Some very interesting relations between the zeros of polynomials,
their derivatives and majorization are studied in a recent paper \cite{BoS} by J.
Borcea and B. Shapiro. They formulate many interesting open problems.


\begin{thebibliography}{123}

\bibitem{Ach} N. I. Achiezer, The classical moment problem, Oliver
and Boyd, Edinburgh, 1965.

\bibitem{Az} A. Aziz, N. A. Rather,
On an inequality of S. Bernstein and the Gauss-Lucas theorem.
Rassias, Themistocles M. (ed.) et al., Analytic and geometric
inequalities and applications. Dordrecht: Kluwer Academic
Publishers. Math. Appl., Dordr. 478, 29-35 (1999)


\bibitem{Berg} T. Bergkvist and H. Rullgard, On polynomial
eigenfunctions for a class of differential operators, Math. Res. Let. v. 9,
p. 153-171, 2002.

\bibitem{BoS} J. Borcea and B. Shapiro, Hiperbolic polynomials and
spectral order, preprint, math. CA/0304145

\bibitem{Cra} T. Craven, G. Csordas, The Gauss-Lucas theorem and
Jensen polynomials. Trans. Amer. Math. Soc. 278 (1983), no. 1, 415--429.

\bibitem{debr1} N. G. de Bruijn and T. A. Springer, On the zeros of a
polynomial and of its derivative II, Indagationes Math. 9, 264-270
(1947).

\bibitem{debr2} M. G. de Bruin, K. G. Ivanov, A. Sharma,
A conjecture of Schoenberg, J. Inequal. Appl. 4, No. 3 (1999), 183-213.

\bibitem{Dim} D. Dimitrov, A refinement of the Gauss-Lucas theorem,
Proc. Am. Math. Soc. 126, No.7 (1998), 2065-2070.

\bibitem{Don} W. F. Donoghue, Monotone Matrix functions and analytic
continuation. Springer, 1974.

\bibitem{FanKy} Ky Fan and G. Pall, Imbedding conditions for Hermitian
and normal matrices, Canad. J. Math. 9 (1957), 298-304.

\bibitem{FiHol} P. Fischer, J. A. R. Holbrook, Balayage defined by
the nonnegative convex functions,
Proc. AMS 79 (1980), 445-448

\bibitem{GeSim} F. Gesztesy and B. Simon,
$m$-functions and inverse spectral analysis for finite and
semi-infinite Jacobi matrices,
J. Anal. Math. 73, 267-297 (1997)

\bibitem{HLP} G. H. Hardy, J. E. Littlwood, G. Polya, Inequalites.
Cambridge, 1988.

\bibitem{Herm} L. H\"ormander, Notions of convexity. Birkh\"auser, 1994.

\bibitem{Hoch} H. Hochstadt, On the construction of a Jacobi
matrix from spectral data, Lin. Algebra and Appl., 8 (1974),
435-446.

\bibitem{HoJ} R. A. Horn and C. R. Johnson, Matrix Analysis.
Cambridge, U.K.: Cambridge Univ. Press, 1993.

\bibitem{Mal1} M. M. Malamud,
On the formula of generalized resolvents of a nondensely
defined Hermitian operator,
Ukr. Matem. Zhurn. vol. 44, N12 (1992), 1658-1688 (in Russian)
(translation in Ukr. Math. J. v.44 (1992), 1522-1547).

\bibitem{Mal2} S. M. Malamud,
Operator inequalities, converse to the Jensen inequality,
Mathematical Notes, v.69, No 4 (2001), 633-637.

\bibitem{Mal3} S. M. Malamud,
A converse to the Jensen inequality, its matrix extensions and
inequalities for minors and eigenvalues,
Linear Algebra and Applications, v. 322, (2001), 19-41.

\bibitem{Mal4} S. M. Malamud, An inverse spectral problem for normal matrices and a generalization of the Gauss-Lucas theorem, math.CV/0304158

\bibitem{Mal5} S. M. Malamud,
Analog of the Poincare separation theorem for normal
matrices and the Gauss-Lukas theorem,
Funct. Anal. Appl., v.37, No 3 (2003), 72-76.

\bibitem{Mar} A. S. Markus, Eigenvalues and singular values of the sum
and product of linear operators,
Russian Math. Surveys 19 (1964), 91-120.

\bibitem{MarM} M. Marcus and H. Minc, A survey on matrix theory and matrix
inequalities. Allyn and Bacon, 1964.

\bibitem{MaOL} A. W. Marshall and I. Olkin, Inequalities: Theory of
majorization and its applications. Acad. Press, 1979.

\bibitem{MS} G. Mason and B. Shapiro, A note on polynomial
eigenfunctions of a hypergeometric type operator, Experimental
mathematics, 10, 609-618.

\bibitem{Mil} M.J. Miller, Maximal polynomials and the Ilieff-Sendov conjecture, Trans. Amer. Math. Soc. 321 (1990), 285--303.

\bibitem{Paw} P Pawlowski,
On the zeros of a polynomial and its derivatives,
Trans. Amer. Math. Soc. 350 (1998), no. 11, 4461--4472.

\bibitem{PoSe} G. Polya and G. Szego, Problems and theorems in
analysis, vol. II, Springer, 1976.

\bibitem{ShaSha} B. Shapiro, M. Shapiro, This strange and misterious Rolle's Theorem, peprint, math.CA/0302215

\bibitem{Shap} Harold S. Shapiro, Spectral aspects of a class of
differential operators, Operator Theory Adv. Appl., 132, p.
361-385. Birkhäuser, Basel, 2002.

\bibitem{Schm} G. Schmeisser, The conjectures of Sendov and Smale,
Approx. Theory: A volume dedicated to Blagovest Sendov, DARBA,
Sofia, 2002, 353-369.

\bibitem{Sher} S. Sherman, On a theorem of Hardy, Littlewood,
Polya, and Blackwell,
Proc. Nat. Acad. Sci. USA 37 (1951), 826-831.

\bibitem{Schoen}
I.J. Schoenberg,
A conjectured analogue of Rolle's theorem for polynomials with
real or complex coefficients,
Am. Math. Mon. 93, 8-13 (1986).

\end{thebibliography}
\end{document}